\documentclass[fleqn,11pt]{wlscirep}
\usepackage[utf8]{inputenc}
\usepackage[T1]{fontenc}
\usepackage{subfigure}
\usepackage{graphicx}
\usepackage{epstopdf}
\usepackage{amsfonts}
\usepackage{csquotes}
\usepackage{amsmath}
\usepackage{amssymb}
\usepackage{graphicx}
\usepackage{subfigure}
\usepackage{tikz}
\usepackage{fancyhdr}
\usepackage{amsmath}
\usepackage{multirow}
\usepackage{color}
\usepackage{float}
\usepackage{array}
\usepackage[english]{babel}
\usepackage{algorithm,algpseudocode}

\makeatletter
\def\BState{\State\hskip-\ALG@thistlm}
\makeatother

\newcommand{\lowertau}[1]{\check{\tau}_{#1}}
\newcommand{\uppertau}[1]{\hat{\tau}_{#1}}
\newcommand{\taustar}[1]{\tau^{*}_{#1}}

\title{Fast generation of stability charts for time-delay systems using continuation of characteristic roots}

\author[1]{Surya Samukham}
\author[2,*]{Thomas K. Uchida}
\author[1]{C. P. Vyasarayani}
\affil[1]{Department of Mechanical and Aerospace Engineering, Indian Institute of Technology Hyderabad, Kandi, Sangareddy 502285, Telangana, India.}
\affil[2]{Department of Mechanical Engineering, University of Ottawa, 161 Louis-Pasteur, Ottawa, Ontario, K1N 6N5, Canada.}
\affil[*]{corresponding: tuchida@uottawa.ca}

\begin{abstract}
Many dynamic processes involve time delays, thus their dynamics are governed by delay differential equations (DDEs).
Studying the stability of dynamic systems is critical, but analyzing the stability of time-delay systems is challenging because DDEs are infinite-dimensional.
We propose a new approach to quickly generate stability charts for DDEs using continuation of characteristic roots (CCR).
In our CCR method, the roots of the characteristic equation of a DDE are written as implicit functions of the parameters of interest, and the continuation equations are derived in the form of ordinary differential equations (ODEs).
Numerical continuation is then employed to determine the characteristic roots at all points in a parametric space; the stability of the original DDE can then be easily determined.
A key advantage of the proposed method is that a system of linearly independent ODEs is solved rather than the typical strategy of solving a large eigenvalue problem at each grid point in the domain.
Thus, the CCR method significantly reduces the computational effort required to determine the stability of DDEs.
As we demonstrate with several examples, the CCR method generates highly accurate stability charts, and does so up to 10 times faster than the Galerkin approximation method.
\end{abstract}

\begin{document}

\flushbottom
\maketitle
\thispagestyle{fancy} 

\section{Introduction}
Many models of dynamic systems involve time delays due to delays in sensing and actuating operations.
Such systems are known as time-delayed systems and their dynamics are governed by delay differential equations (DDEs).
DDEs have been investigated extensively in recent years due to their wide-ranging applications in modeling a large number of natural and control processes~\cite{pekar2018spectrum}.
Some examples include control systems~\cite{sieber2008control}, manufacturing~\cite{kalmar2001subcritical,balachandran2001nonlinear,insperger2003multiple,insperger2004stability,long2010stability,nayfeh2012time}, lasers~\cite{kane2005unlocking}, the delayed feedback control mechanism of human balancing~\cite{foss1996multistability,stepan2000balancing,ahsan2016balance}, traffic flow models~\cite{orosz2004hopf}, biology~\cite{bocharov2000numerical,popovych2017closed}, epilepsy seizure models~\cite{rodrigues2009transitions}, physics~\cite{kantner2015delay,alvarez2017advanced}, and many other engineering applications~\cite{kyrychko2010use}.
Recently, Young et al.~\cite{young2019consequences} studied the consequences of delays and imperfect implementation of isolation in epidemic control using time-delayed dynamic system models.

A critical study for any dynamic system is analyzing its stability.
In stable regions of a parametric space, small perturbations decay over time and the system remains ``well-behaved''; in unstable regions, the dynamics of the system diverge with potentially disastrous consequences.
Determining the stability of a DDE, or the regions of stability in a parametric space, is challenging because DDEs are infinite-dimensional~\cite{stepan1989retarded,olgac2002exact}.
One strategy to determine the stability of a DDE is to compute the locations of its characteristic roots in the complex plane.
The characteristic equation of a DDE is a quasi-polynomial with infinitely many roots; the DDE is stable if, and only if, all the roots lie in the left half of the complex plane.
In the literature, several methods have been proposed to approximate the characteristic roots of DDEs for studying their stability.
Some examples include the semi-discretization method~\cite{insperger2002semi}, D-subdivision methods~\cite{olgac2002exact}, finite difference methods~\cite{sun2009control}, finite element methods~\cite{mann2010stability}, mapping-based algorithms for large-scale computation of quasi-polynomial roots~\cite{vyhlidal2009mapping}, and Galerkin approximations~\cite{wahi2005galerkin,vyasarayani2012galerkin}.
The Lambert W function is another powerful technique to determine the stability of DDEs, however it can be used only when a single delay is present~\cite{yi2010timedelay}.
In most of the aforementioned methods, the characteristic roots of the DDE are evaluated by solving an eigenvalue problem.
Therefore, to find regions of stability in a parametric space, the region must first be discretized into a finite grid of sufficient density, and then an eigenvalue problem must be solved at each grid point.
This approach requires substantial computational effort and is not an ideal strategy to determine stability regions or boundaries with high accuracy.

Methods have also been developed to determine the stability of a DDE without calculating its characteristic roots.
For example, the direct numerical integration of a DDE provides its time response and therefore reveals its stability.
However, to determine the regions of stability in a parametric space using this method, it would be necessary to analyze the time response of the DDE at all points in the spectrum.
Analytical stability boundaries can be obtained by tracking all the critical curves on which at least one pair of purely imaginary roots exists.
However, this method does not provide any information about the stable and unstable regions in the spectrum.
Also, it cannot be guaranteed that the critical curves always represent the stability boundary: it may happen that a pair of characteristic roots lies on the imaginary axis while another pair lies in the right half of the complex plane, in which case the system is unstable.
Cluster treatment of characteristic roots~\cite{olgac2004cluster} can be used to generate exact stability charts for DDEs; however, this strategy does not provide any information about the characteristic roots or their locations.
More recently, Che et al.~\cite{che2019multi} proposed a multi-fidelity model for identifying the stability boundary in time-delayed systems.
In this approach, the stability boundaries are identified accurately by refining the mesh at the critical regions.
While the computational effort required for the approach of Che et al.\ is small relative to many other methods, a large number of eigenvalue problems must still be solved to determine the stability boundaries.

The methods discussed above suggest that determining the stability of a DDE, or the regions of stability in a parametric space, is a computationally expensive task.
In this work, we have developed a continuation of characteristic roots (CCR) method to determine the characteristic roots and thus the stability regions of DDEs with relatively low computational cost.
In the CCR method, we first write the characteristic roots as implicit functions of the parameters of interest and derive the continuation equations in the form of ordinary differential equations (ODEs), using the chain rule of differentiation.
Upon solving these ODEs using appropriate initial conditions, we obtain the corresponding roots with respect to the parameters of interest.
Thus, very accurate stability charts are obtained simply by solving systems of linearly independent ODEs rather than solving a large number of eigenvalue problems.

This paper is organized as follows.
In Sec.~\ref{sec.modeling}, we describe the CCR method for determining the roots of the characteristic equation of a DDE and the strategy to determine its regions of stability.
In Sec.~\ref{sec.results}, we provide several examples to demonstrate the efficacy of the proposed CCR method.
We also discuss a scenario in which the method fails and recommend a technique to address this limitation.
Finally, we summarize our findings in Sec.~\ref{sec.conclusions}.

\section{Mathematical Modeling}\label{sec.modeling}
In this section, we describe the mathematical procedure for applying the CCR method by considering a second-order DDE of the following form:
\begin{equation}
    \ddot{x}(t) + a_1 \dot{x}(t) + a_2 x(t) + b_1 \dot{x}(t-\tau_1) + b_2 x(t-\tau_2) = 0,
    \label{mm_eq1}
\end{equation}
where $x(t)$ is the system state vector, $\dot{x}(t)$ and $\ddot{x}(t)$ are its first and second derivatives with respect to time, $\{a_1, a_2, b_1, b_2\} \in \mathbb{R}$ are parameters, and time delays $\tau_i \ge 0$ for $i=1,~2$.
Equation~\eqref{mm_eq1} is a DDE if any $\tau_i > 0$; otherwise, it is simply an ODE.
History functions that describe the past states of the system are given as follows:
\begin{flalign}
    x(t) &= \alpha(t),\\
    \dot{x}(t) &= \beta(t), \quad -\bar{\tau} \le t \le 0,
    \label{mm_eq2}
\end{flalign}
where $\bar{\tau} \triangleq \text{max}(\tau_1, \tau_2)$.
The characteristic equation of the DDE is obtained by substituting $x = e^{\lambda t}$ into Eq.~\eqref{mm_eq1}:
\begin{equation}
    D(\lambda) \triangleq \lambda^2 + a_1 \lambda + a_2 + b_1 \lambda e^{-\lambda \tau_1} + b_2 e^{-\lambda \tau_2} = 0.
    \label{mm_eq3}
\end{equation}
To determine the characteristic roots ($\lambda$) corresponding to the first time delay ($\tau_1$) in Eq.~\eqref{mm_eq3}, we write $\lambda$ as an implicit function of $\tau_1$ and, from the chain rule of differentiation, we have the following:
\begin{gather}
    dD(\lambda, \tau_1) \equiv \frac{\partial D}{\partial \lambda} d\lambda + \frac{\partial D}{\partial \tau_1} d\tau_1 = 0,\\
    \frac{d\lambda}{d\tau_1} = - \frac{\partial D}{\partial \tau_1} \bigg/ \frac{\partial D}{\partial \lambda} \implies \frac{b_1 \lambda^2 e^{-\lambda \tau_1}}{2\lambda + a_1 + b_1 \left( 1 - \lambda \tau_1 \right) e^{-\lambda \tau_1} - b_2 \tau_2 e^{-\lambda \tau_2}}.
    \label{mm_eq4}
\end{gather}
Similarly, to determine $\lambda$ corresponding to $\tau_2$, we write $\lambda$ as an implicit function of $\tau_2$ and proceed as above:
\begin{gather}
    dD(\lambda, \tau_2) \equiv \frac{\partial D}{\partial \lambda} d\lambda + \frac{\partial D}{\partial \tau_2} d\tau_2 = 0,\\
    \frac{d\lambda}{d\tau_2} = - \frac{\partial D}{\partial \tau_2} \bigg/ \frac{\partial D}{\partial \lambda} \implies \frac{b_2 \lambda e^{-\lambda \tau_2}}{2\lambda + a_1 + b_1 \left( 1 - \lambda \tau_1 \right) e^{-\lambda \tau_1} - b_2 \tau_2 e^{-\lambda \tau_2}}.
    \label{mm_eq5}
\end{gather}
Upon solving the ODEs given by Eqs.~\eqref{mm_eq4} and \eqref{mm_eq5}, we obtain the roots of the characteristic equation (Eq.~\eqref{mm_eq3}) corresponding to delays $\tau_1$ and $\tau_2$, respectively.
The initial conditions (roots) to solve the ODEs (Eqs.~\eqref{mm_eq4} and \eqref{mm_eq5}) can be obtained using any of various existing methods~\cite{insperger2002semi,vyasarayani2012galerkin}; in this work, we use Galerkin approximation~\cite{wahi2005galerkin,vyasarayani2012galerkin} to determine the initial roots.
To compute multiple roots simultaneously, Eqs.~\eqref{mm_eq4} and \eqref{mm_eq5} are written as a system of linearly independent ODEs:
\begin{equation}
    \frac{d\boldsymbol{\lambda}}{d \tau_i} = - \boldsymbol{J}^{-1} \frac{\partial \mathbf{D}_i}{\partial \tau_i}, \quad i = 1, 2,
    \label{mm_eq6}
\end{equation}
where $\boldsymbol{\lambda} = \left[ \lambda_1, \lambda_2, \ldots, \lambda_N \right]^{\operatorname{T}}$ is a vector of characteristic roots, $\boldsymbol{J}$ is a Jacobian matrix given by $\boldsymbol{J} = \operatorname{diag} \left( \frac{\partial D}{\partial \lambda_1}, \frac{\partial D}{\partial \lambda_2}, \ldots, \frac{\partial D}{\partial \lambda_N} \right)$, and $\mathbf{D}_i$ is a diagonal matrix given by $\mathbf{D}_i = \operatorname{diag} \left( D(\lambda_1,\tau_i), D(\lambda_2,\tau_i), \ldots, D(\lambda_N,\tau_i) \right)$.
The system of ODEs in Eq.~\eqref{mm_eq6} is solved, using the roots obtained from the Galerkin approximation method as initial conditions, to determine the corresponding roots with respect to parameter $\tau_i$.

Suppose we wish to determine the stability regions of the DDE (Eq.~\eqref{mm_eq1}) in the parametric space of $\tau_1 \in \left[ \lowertau{1}, \uppertau{1} \right]$ and $\tau_2 \in \left[ \lowertau{2}, \uppertau{2} \right]$.
We first use the Galerkin approach to evaluate the $N$ rightmost characteristic roots for Eq.~\eqref{mm_eq3} at any point $\left( \taustar{1}, \taustar{2} \right)$ in the parametric space, where $\lowertau{1} \leq \taustar{1} \leq \uppertau{1}$ and $\lowertau{2} \leq \taustar{2} \leq \uppertau{2}$.
The obtained roots are then used as initial conditions to solve the system of ODEs (Eq.~\eqref{mm_eq6}) over the domains $\tau_1 \in \left[ \lowertau{1}, \taustar{1} \right]$ and $\tau_1 \in \left[ \taustar{1}, \uppertau{1} \right]$, holding $\tau_2 = \taustar{2}$ constant.
The solution is then evaluated at specified grid points where $\tau_1 \in \left[ \lowertau{1}, \uppertau{1} \right]$ and $\tau_2 = \taustar{2}$.
Note that, in the domain of integration $\tau_1 \in \left[ \lowertau{1}, \taustar{1} \right]$, we begin at $\tau_1 = \taustar{1}$ and solve for decreasing $\tau_1$.
Upon completion of this stage, we have obtained through numerical continuation the corresponding $N$ characteristic roots at each point in the domain $\tau_1 \in \left[ \lowertau{1}, \uppertau{1} \right]$ and $\tau_2 = \taustar{2}$.
Next, we use each of these solutions as initial conditions to solve the system of ODEs (Eq.~\eqref{mm_eq6}) along the $\tau_2$ dimension---that is, over the domains $\tau_2 \in \left[ \lowertau{2}, \taustar{2} \right]$ and $\tau_2 \in \left[ \taustar{2}, \uppertau{2} \right]$---while holding $\tau_1$ constant in each integration.
We repeat for each solution along $\tau_1 \in \left[ \lowertau{1}, \uppertau{1} \right]$ computed earlier.
Upon completion of this stage, we have obtained the corresponding characteristic roots at all points in the parametric space $\left[ \lowertau{1}, \uppertau{1} \right] \times \left[ \lowertau{2}, \uppertau{2} \right]$.
The stability charts for the original DDE system (Eq.~\eqref{mm_eq1}) can then be generated simply by determining the location in the complex plane of the rightmost characteristic root at each grid point in the parametric space.
The CCR method has been summarized in Algorithm~1.

\begin{algorithm}
    \caption{The CCR method for fast generation of stability charts for DDEs}
    \begin{algorithmic}
        \Statex
        \State{\textbf{Given:}} A second-order DDE of the form given in Eq.~\eqref{mm_eq1} and a parametric space $\left[ \lowertau{1}, \uppertau{1} \right] \times \left[ \lowertau{2}, \uppertau{2} \right]$.
        \State{\textbf{Find:}} The stability chart for the DDE over the specified parametric space.
        \State{$\boldsymbol{\lambda}_{\text{IC}}$ $\gets$ $N$ rightmost characteristic roots at any point $\left( \taustar{1}, \taustar{2} \right)$ in the parametric space, obtained using the Galerkin approach.}
        \State{$\boldsymbol{\lambda}_{\tau_1}$ $\gets$ solve $\frac{d \boldsymbol{\lambda}}{d \tau_1} = - \boldsymbol{J}^{-1}\frac{\partial \mathbf{D}_1}{\partial \tau_1}$ over the domains $\tau_1 \in \left[ \lowertau{1}, \taustar{1} \right]$ and $\tau_1 \in \left[ \taustar{1}, \uppertau{1} \right]$, using $\boldsymbol{\lambda}_{\text{IC}}$ as initial conditions and holding $\tau_2 = \taustar{2}$ constant, and evaluate the solution at $N_{\tau_1}$ grid points.}
        \For{$i$ from $1$ to $N_{\tau_1}$}
            \State{$\boldsymbol{\lambda}_{\tau_2}^i$ $\gets$ solve $\frac{d \boldsymbol{\lambda}}{d \tau_2} = - \boldsymbol{J}^{-1}\frac{\partial \mathbf{D}_2}{\partial \tau_2}$ over the domains $\tau_2 \in \left[ \lowertau{2}, \taustar{2} \right]$ and $\tau_2 \in \left[ \taustar{2}, \uppertau{2} \right]$, using $\boldsymbol{\lambda}_{\tau_1}^i$ as initial conditions and holding $\tau_1 = \tau_1^i$ constant, and evaluate the solution at $N_{\tau_2}$ grid points.}
        \EndFor\\
        \# Check for stability
        \For{$i$ from $1$ to $N_{\tau_1}$}
            \For{$j$ from $1$ to $N_{\tau_2}$}
                \If{$\text{max}\left( \text{Re}\left\{ \boldsymbol{\lambda}_{\tau_2}^{i,j} \right\} \right) \leq 0$}
                    \State{Stable at the point $\left( \tau_1^i, \tau_2^j \right)$.}
                \Else
                    \State{Unstable at the point $\left( \tau_1^i, \tau_2^j \right)$.}
                \EndIf
            \EndFor
        \EndFor
    \end{algorithmic}
\end{algorithm}

\section{Results}\label{sec.results}
In this section, we generate the stability charts for three DDEs using the proposed CCR method and compare them with the stability charts obtained using the Galerkin approach.
We discuss the root-crossing phenomenon and demonstrate why several characteristic roots must be continued to obtain accurate stability charts; simply continuing the rightmost root is inadequate.
Finally, we present a scenario in which continuation fails and provide a strategy to address this limitation.

\subsection{Example 1}
We first consider the following first-order DDE with five delays:
\begin{equation}
    \dot{x} + ax + \sum_{i=1}^{5} b_i x(t-\tau_i) = 0.
    \label{fo_eq1}
\end{equation}
Upon substituting $x = e^{\lambda t}$ into Eq.~\eqref{fo_eq1}, we obtain the following characteristic equation for the DDE:
\begin{equation}
    D(\lambda) \equiv \lambda + a + \sum_{i=1}^{5} b_i e^{-\lambda \tau_i} = 0.
    \label{fo_eq2}
\end{equation}
By considering $\lambda$ as an implicit function of $\tau_1$ and following a similar mathematical approach as described in Sec.~\ref{sec.modeling}, we arrive at the following ODE:
\begin{equation}
    \frac{d \lambda}{d \tau_1} = -\frac{\partial D}{\partial \tau_1} \bigg/ \frac{\partial D}{\partial \lambda} \implies \frac{b_1 \lambda e^{-\lambda \tau_1}}{1 - \sum_{i=1}^{5} b_i \tau_i e^{-\lambda \tau_i}}.
    \label{fo_eq3}
\end{equation}
Upon solving the above ODE (Eq.~\eqref{fo_eq3}), using the rightmost characteristic roots obtained from the Galerkin approach as initial conditions, we obtain the roots corresponding to $\tau_1$.

\subsubsection*{Root-crossing}\label{rootcrossing}
\begin{figure}[htp]
    \begin{center}
    \includegraphics[width=0.45\textwidth]{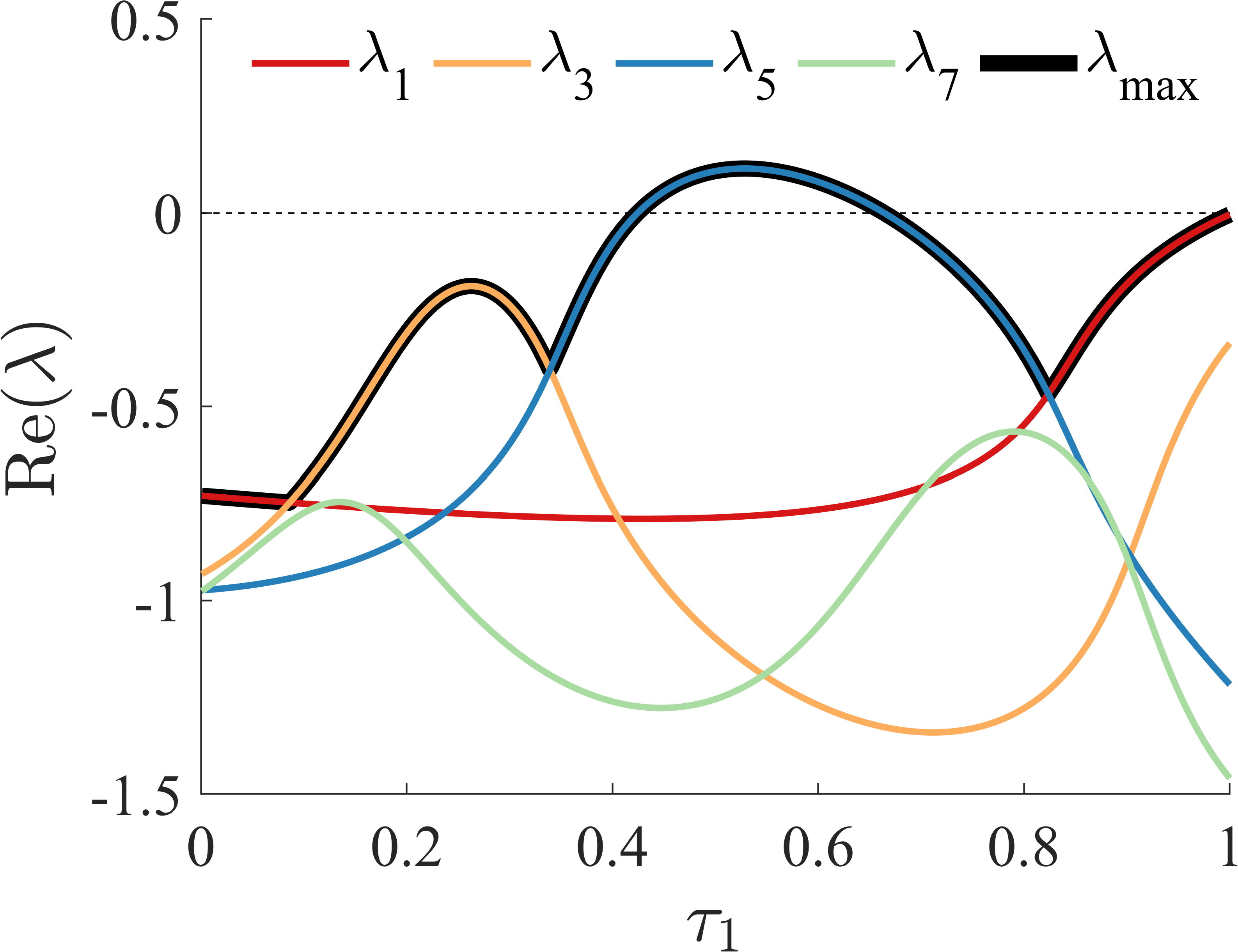}
    \end{center}
    \caption{Crossing of the characteristic roots of the DDE given by Eq.~\eqref{fo_eq1}.}
    \label{fig_rootcross}
\end{figure}

In this section, we explore the accuracy of the CCR method at determining the rightmost characteristic root for the DDE.
We first determine the 8 rightmost roots for the characteristic equation (Eq.~\eqref{fo_eq2}) using the Galerkin approach.
We then use these roots as initial conditions to solve the ODE (Eq.~\eqref{fo_eq3}) and determine $\lambda$ in the domain $\tau_1 \in \left[ 0.001, 1 \right]$.
The parameters used for this analysis are as follows: $a=1$, $b_1=3$, $b_2=2.8$, $b_3=0.6$, $b_4=0.8$, $b_5=1$, $\taustar{1}=0.001$, $\tau_2=0.25$, $\tau_3=1$, $\tau_4=1.5$, and $\tau_5=2$.
The real part of each characteristic root obtained from the continuation method is shown in Fig.~\ref{fig_rootcross}.
Because the characteristic roots appear as complex conjugates, only the odd-numbered roots are shown.
The dominant (rightmost) characteristic root corresponding to $\tau_1$ obtained using the Galerkin approach ($\lambda_{\text{max}}$) has also been shown in Fig.~\ref{fig_rootcross}.
In the Galerkin approach, $\lambda_{\text{max}}$ was obtained by discretizing the domain $\tau_1 \in \left[ 0.001, 1 \right]$ into 200 grid points and solving an eigenvalue problem at each point.

This example demonstrates that the dominant characteristic root at the initial point may not be the dominant root throughout the domain.
In Fig.~\ref{fig_rootcross}, the dominant root at the initial point $\tau_1 = 0.001$ (i.e., $\lambda_1$) is dominant only in the domain $\tau_1 \in \left[ 0.001, 0.08 \right]$.
At $\tau_1 = 0.08$, the third root ($\lambda_3$) crosses $\lambda_1$ and is dominant in the interval $\tau_1 = \left[ 0.08, 0.33 \right]$; the fifth root ($\lambda_5$) is dominant in the interval $\tau_1 \in \left[ 0.33, 0.82 \right]$ and, finally, the first root ($\lambda_1$) is again dominant in the interval $\tau_1 \in \left[ 0.82, 1 \right]$.
Figure~\ref{fig_rootcross} clearly illustrates the root-crossing phenomenon in the characteristic roots of the DDE and proves that the rightmost root at one point may not remain dominant throughout the domain.
The location of the rightmost root determines whether the system is stable and, thus, is critical information for generating the stability charts for a DDE.
Therefore, we observe that it is not sufficient to continue only the rightmost root; we instead track the $N$ rightmost roots and increase $N$ until convergence is achieved.

\subsubsection*{Stability chart}\label{stabilitychart}
We now use the CCR method to determine the regions of stability for the DDE (Eq.~\eqref{fo_eq1}) in the parametric space $\tau_1 \in \left[0.001, 1 \right]$ and $\tau_2 \in \left[ 0.001, 1 \right]$.
To generate the stability chart, we first use Galerkin approximation to obtain the 25 rightmost characteristic roots for the DDE at the initial point $\left( \tau_1, \tau_2 \right) = \left( 0.001, 0.001 \right)$.
The other parameters are as follows: $a=1$, $b_1=3$, $b_2=2.8$, $b_3=0.6$, $b_4=0.8$, $b_5=1$, $\tau_3=1$, $\tau_4=1.5$, and $\tau_5=2$.
To determine the roots at the initial point with high accuracy, we use $N_G = 200$ modes in the Galerkin approximation.
We use these roots as initial conditions to solve the ODE (Eq.~\eqref{fo_eq3}) and determine the characteristic roots that correspond to $\tau_1$ in the interval $\tau_1 \in \left[ 0.001, 1 \right]$.
The integration is performed in \textsc{Matlab} using the ``ode45'' explicit integrator with absolute and relative tolerances of $10^{-12}$.
The solution of the ODEs is then evaluated at 2000 equidistant points in the interval $\tau_1 \in \left[ 0.001, 1 \right]$ to obtain the corresponding 25 roots for Eq.~\eqref{fo_eq2}.
Next, we use the obtained roots from Eq.~\eqref{fo_eq3} as initial conditions and continue the roots with respect to $\tau_2$ using the following equation:
\begin{equation}
    \frac{d \lambda}{d \tau_2} = -\frac{\partial D}{\partial \tau_2} \bigg/ \frac{\partial D}{\partial \lambda} \implies \frac{b_2 \lambda e^{-\lambda \tau_2}}{1 - \sum_{i=1}^{5} b_i \tau_i e^{-\lambda \tau_i}}.
    \label{fo_eq4}
\end{equation}
Equation~\eqref{fo_eq4} is solved and evaluated at 2000 equidistant grid points in the interval $\tau_2 \in \left[ 0.001, 1 \right]$ for each point along $\tau_1 \in \left[ 0.001, 1 \right]$.
Following this integration step, the corresponding characteristic roots of the DDE (Eq.~\eqref{fo_eq1}) at all grid points in the parametric space $\tau_1 \in \left[ 0.001, 1 \right]$ and $\tau_2 \in \left[ 0.001, 1 \right]$ have been determined.
Finally, we obtain the stability regions of the DDE by analyzing the location of the characteristic roots in the complex plane at each point in the parametric space.

\begin{figure}[htp]
    \begin{center}
    \includegraphics[width=0.9\textwidth]{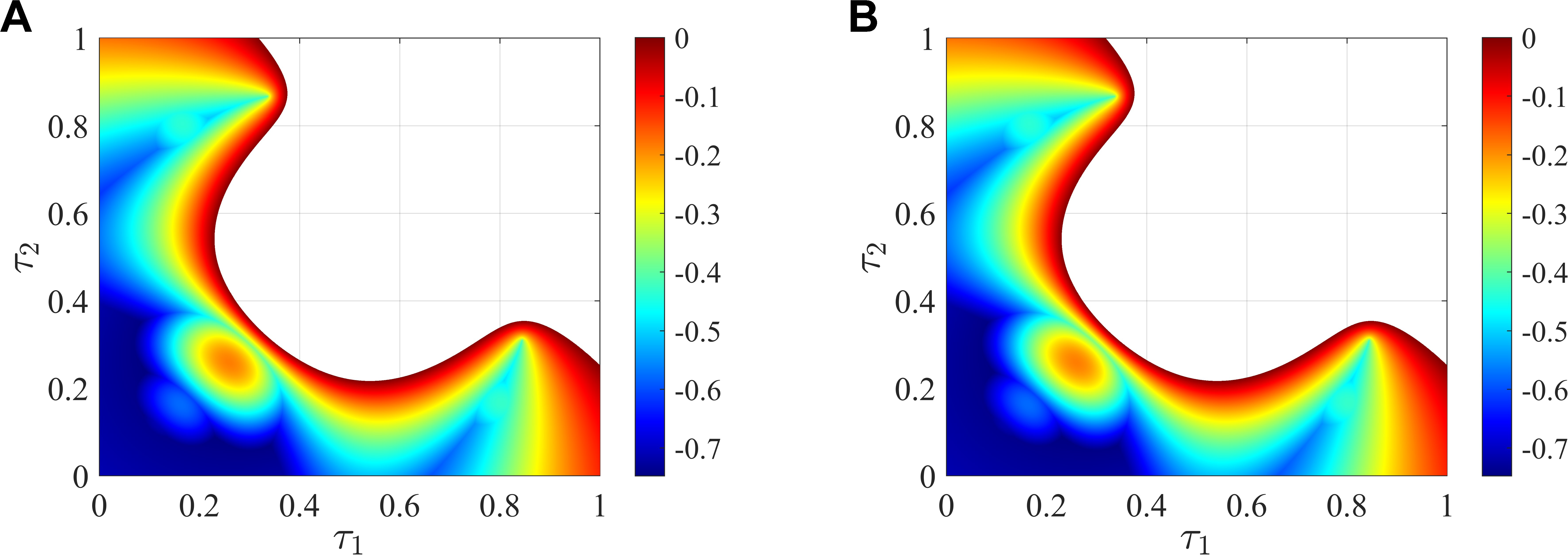}
    \end{center}
    \caption{Stability chart for the first-order DDE with five delays (Eq.~\eqref{fo_eq1}) obtained using (A) the CCR method and (B) the Galerkin approach.}
    \label{ex1fig1}
\end{figure}

The stability chart obtained using the CCR method is shown in Fig.~\ref{ex1fig1}(A).
The color contours in the figure represent the maximum damping present in the system (i.e., the real part of the dominant characteristic root).
To verify the results obtained from the CCR method, we also present the results obtained using Galerkin approximation only, shown in Fig.~\ref{ex1fig1}(B).
In the Galerkin-only approach, we discretize the parametric space into a $2000 \times 2000$ grid and solve an eigenvalue problem at each grid point using $N_G = 25$ to determine the dominant characteristic root.
The results presented in Fig.~\ref{ex1fig1} clearly demonstrate the correctness of the results obtained using the CCR method.

\subsection{Example 2}
\begin{figure}[htp]
    \begin{center}
    \includegraphics[width=0.9\textwidth]{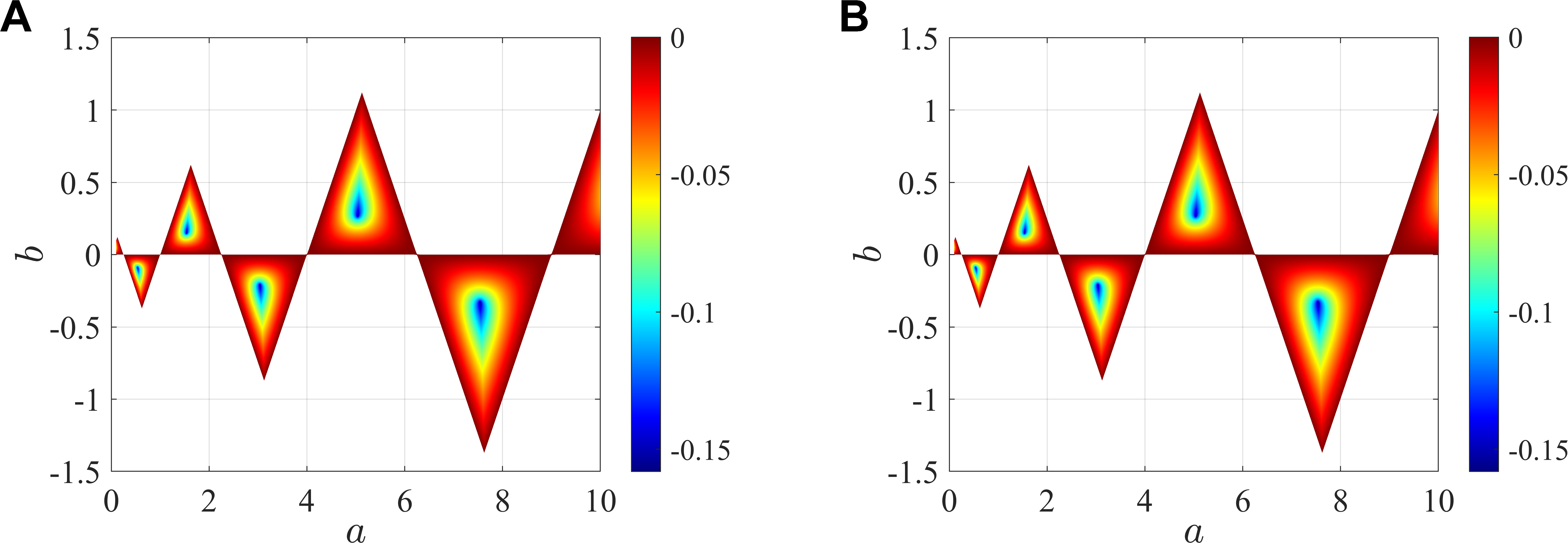}
    \end{center}
    \caption{Stability chart for the second-order DDE with a single delay (Eq.~\eqref{SO_eq1}) obtained using (A) the CCR method and (B) the Galerkin approach.}
    \label{ex2fig1}
\end{figure}

In this example, we consider the following second-order DDE with a single delay:
\begin{equation}
    \ddot{x}(t) + a x(t) - b x(t-\tau) = 0,
    \label{SO_eq1}
\end{equation}
which has the following characteristic equation:
\begin{equation}
    D(\lambda) \equiv \lambda^2 + a - b e^{-\lambda \tau} = 0.
    \label{SO_eq2}
\end{equation}
We determine the stability regions of the DDE (Eq.~\eqref{SO_eq1}) in the parametric space $a \in \left[ 0.01, 10 \right]$ and $b \in \left[ -1.5, 1.5 \right]$.
We write $\lambda$ as an implicit function of $a$ and $b$ separately, and write the continuation differential equations using the chain rule of differentiation as follows:
\begin{flalign}
    \frac{d \lambda}{d a} &= -\frac{\partial D}{\partial a} \bigg/ \frac{\partial D}{\partial \lambda} \implies -\frac{1}{2\lambda + b \tau e^{-\lambda \tau}},\label{SO_eq4a}\\
    \frac{d \lambda}{d b} &= -\frac{\partial D}{\partial b} \bigg/ \frac{\partial D}{\partial \lambda} \implies \frac{e^{-\lambda \tau}}{2\lambda + b \tau e^{-\lambda \tau}}.\label{SO_eq4b}
\end{flalign}
As in the previous example, we begin by determining the 25 rightmost characteristic roots for Eq.~\eqref{SO_eq2} using the Galerkin approach, in this case using the initial point $\left(a, b \right) = \left(0.01, -1.5 \right)$; the time delay parameter in Eq.~\eqref{SO_eq1} is set to $\tau = 2\pi$.
We then use the roots obtained from the Galerkin approach as initial conditions to solve the ODE in Eq.~\eqref{SO_eq4a} over the domain $a \in \left[ 0.01, 10 \right]$ with $b = -1.5$ held constant.
The roots obtained from Eq.~\eqref{SO_eq4a} are then used as initial conditions to solve Eq.~\eqref{SO_eq4b} along $b \in \left[ -1.5, 1.5 \right]$ for each point in the domain $a \in \left[ 0.01, 10 \right]$.
The stability regions thus obtained from the CCR method are shown in Fig.~\ref{ex2fig1}(A); the results obtained using the Galerkin approach are shown in Fig.~\ref{ex2fig1}(B) for verification.
A grid size of $2000 \times 2000$ was used for both methods.
The results presented in Fig.~\ref{ex2fig1} again demonstrate the correctness of the results obtained using the CCR method.

\subsection{Example 3}
\begin{figure}[htp]
    \begin{center}
    \includegraphics[width=0.9\textwidth]{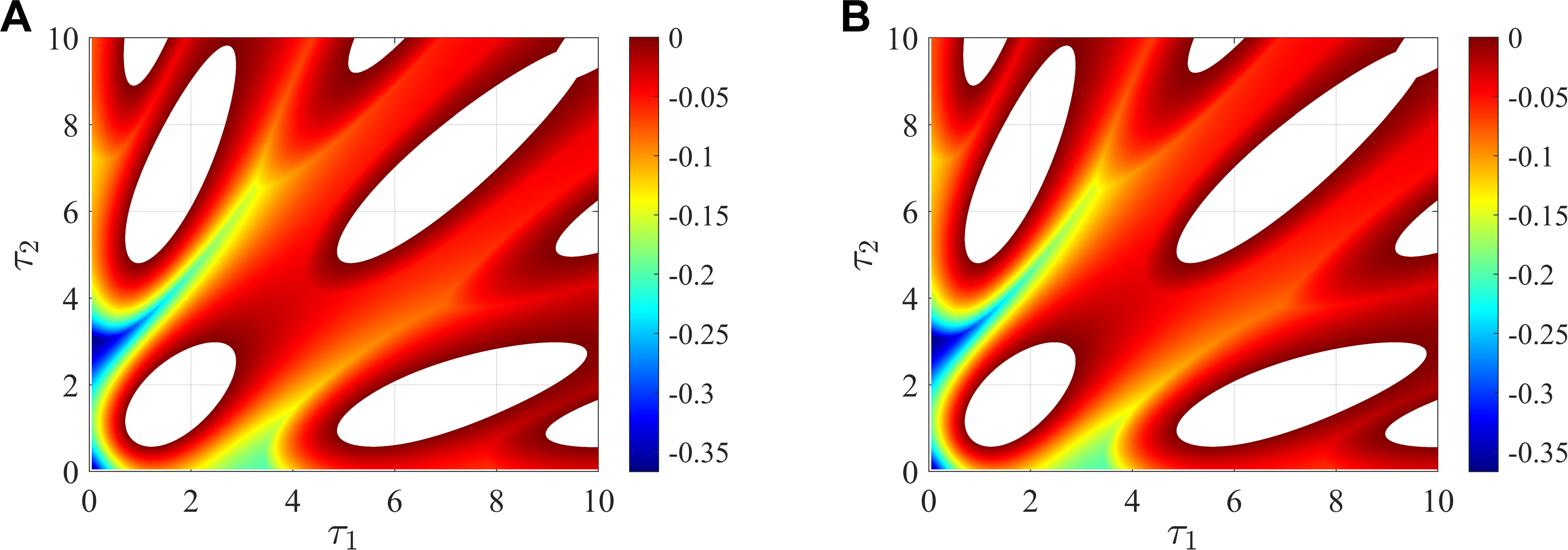}
    \end{center}
    \caption{Stability chart for the second-order DDE with two delays (Eq.~\eqref{SO_2delaywd_eq1}) using parameter set 1, obtained using (A) the CCR method and (B) the Galerkin approach.}
    \label{ex3fig1}
\end{figure}

\begin{figure}[htp]
    \begin{center}
    \includegraphics[width=0.9\textwidth]{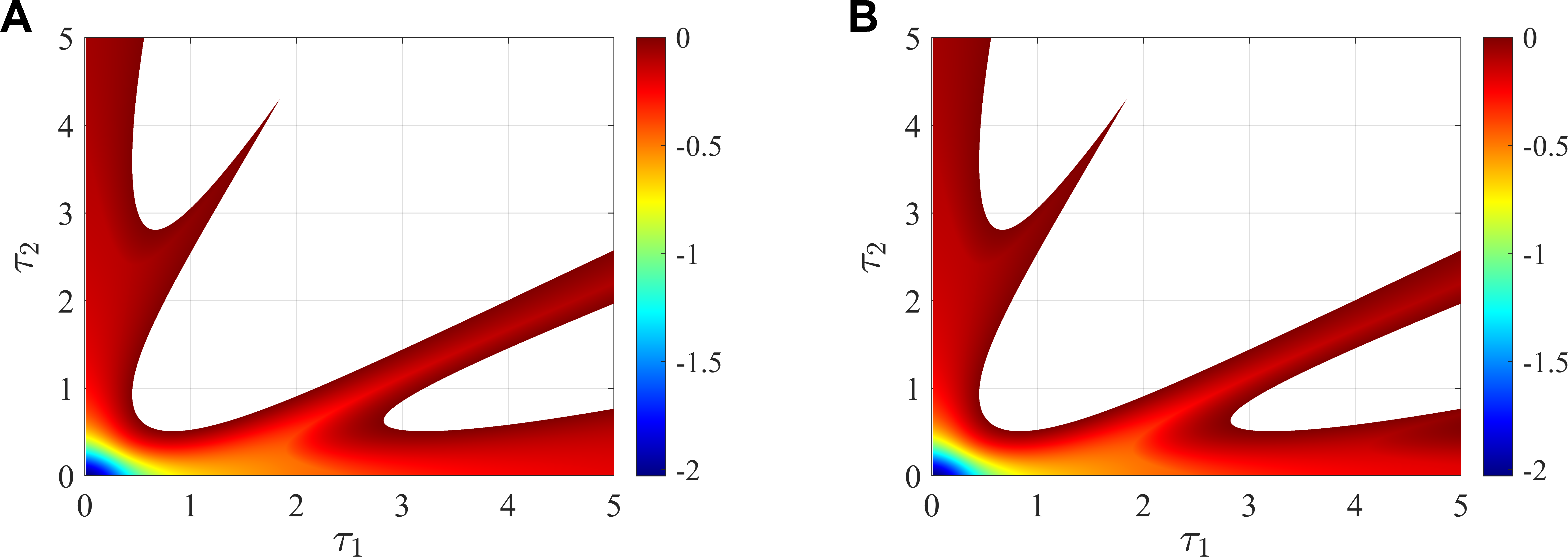}
    \end{center}
    \caption{Stability chart for the second-order DDE with two delays (Eq.~\eqref{SO_2delaywd_eq1}) using parameter set 2, obtained using (A) the CCR method and (B) the Galerkin approach.}
    \label{ex3fig2}
\end{figure}

We now consider the following second-order DDE with two delays:
\begin{equation}
    \ddot{x}(t) + a_1 \dot{x}(t) + a_2 x(t) + b_1 \dot{x}(t-\tau_1) + b_2 x(t-\tau_1) + b_3 \dot{x}(t-\tau_2) + b_4 x(t-\tau_2) = 0,
    \label{SO_2delaywd_eq1}
\end{equation}
which has the following characteristic equation:
\begin{equation}
    D(\lambda) \equiv \lambda^2 + a_1 \lambda + a_2 + b_1 \lambda e^{-\lambda \tau_1} + b_2 e^{-\lambda \tau_1} + b_3 \lambda e^{-\lambda \tau_2} + b_4 e^{-\lambda \tau_2} = 0.
    \label{SO_2delay_eq2}
\end{equation}
To determine the stability regions of the DDE (Eq.~\eqref{SO_2delaywd_eq1}) in the parametric space of $\tau_1$ and $\tau_2$, we write $\lambda$ as an implicit function of $\tau_1$ and $\tau_2$ separately, and derive the continuation differential equations using the chain rule of differentiation:
\begin{flalign}
    \frac{d \lambda}{d \tau_1} &= -\frac{\partial D}{\partial \tau_1} \bigg/ \frac{\partial D}{\partial \lambda} \implies \frac{\left( b_1 \lambda^2 + b_2 \lambda \right) e^{-\lambda \tau_1}}{\Delta},\label{SO_2delay_eq3a}\\
    \frac{d \lambda}{d \tau_2} &= -\frac{\partial D}{\partial \tau_2} \bigg/ \frac{\partial D}{\partial \lambda} \implies \frac{\left( b_3 \lambda^2 + b_4 \lambda \right) e^{-\lambda \tau_2}}{\Delta},\label{SO_2delay_eq3b}
\end{flalign}
where
\begin{equation}
    \Delta = 2\lambda + a_1 + b_1 \left( 1 - \lambda \tau_1 \right) e^{-\lambda \tau_1} - b_2 \tau_1 e^{-\lambda \tau_1} + b_3 \left( 1 - \lambda \tau_2 \right) e^{-\lambda \tau_2} - b_4 \tau_2 e^{-\lambda \tau_2}.
\end{equation}
For the DDE given by Eq.~\eqref{SO_2delaywd_eq1}, we analyze the stability of the system in the parametric space of $\tau_1$ and $\tau_2$ for two sets of parameters:
\begin{itemize}
    \item Set 1: $a_1=0.8$, $a_2=1.9$, $b_1=0$, $b_2=0.8$, $b_3=0$, and $b_4=0.5$
    \item Set 2: $a_1=3$, $a_2=5$, $b_1=0.5$, $b_2=3$, $b_3=0.6$, and $b_4=5.2$
\end{itemize}
The stability analysis for the DDE is performed using the CCR method by solving the ODEs given by Eqs.~\eqref{SO_2delay_eq3a} and \eqref{SO_2delay_eq3b}.
The initial conditions for the ODEs at the initial point $\left( \tau_1, \tau_2 \right) = \left( 0.01, 0.01 \right)$ are obtained using the Galerkin approach.
The stability charts generated using the CCR method for parameter sets 1 and 2 are shown in Figs.~\ref{ex3fig1}(A) and \ref{ex3fig2}(A), respectively; the corresponding stability charts obtained using Galerkin approach are shown in Figs.~\ref{ex3fig1}(B) and \ref{ex3fig2}(B) for comparison.
All stability charts in Figs.~\ref{ex3fig1} and \ref{ex3fig2} are generated over a grid size of $2000 \times 2000$.
Once again, the stability regions found using the CCR method match those found using the Galerkin approximation method.

\subsection{Rank-deficient Jacobian}\label{failcase}
While the above examples demonstrate the efficacy of the proposed CCR method for determining the characteristic roots and stability regions of a DDE, the method has a limitation.
It is possible that the derived continuation ODEs will become discontinuous for certain parameters, in which case the solution of the ODE cannot be determined.
One such case has been encountered for the DDE given by Eq.~\eqref{SO_2delaywd_eq1} using the following parameters (set 3): $a_1=1.5$, $a_2=0.8$, $b_1=2$, $b_2=0.5$, $b_3=1$, and $b_4=1$.
When these parameters are used, $\frac{\partial \lambda}{\partial \tau_i} = 0$ for certain combinations of $\tau_1$, $\tau_2$, and $\lambda$, leading to $\left| \boldsymbol{J} \right| = 0$ and, thus, the differential equations given by Eqs.~\eqref{SO_2delay_eq3a} and \eqref{SO_2delay_eq3b} become discontinuous.
When a \textsc{Matlab} integrator is used to solve such ODEs, it fails to proceed when $\left| \boldsymbol{J} \right| = 0$.
To overcome this limitation, we terminate the integration whenever $\left| \boldsymbol{J} \right| = 0$ and resume integration at the next grid point with a new set of initial conditions evaluated using the Galerkin approach at the corresponding point.

\begin{figure}[htp]
    \begin{center}
    \includegraphics[width=0.9\textwidth]{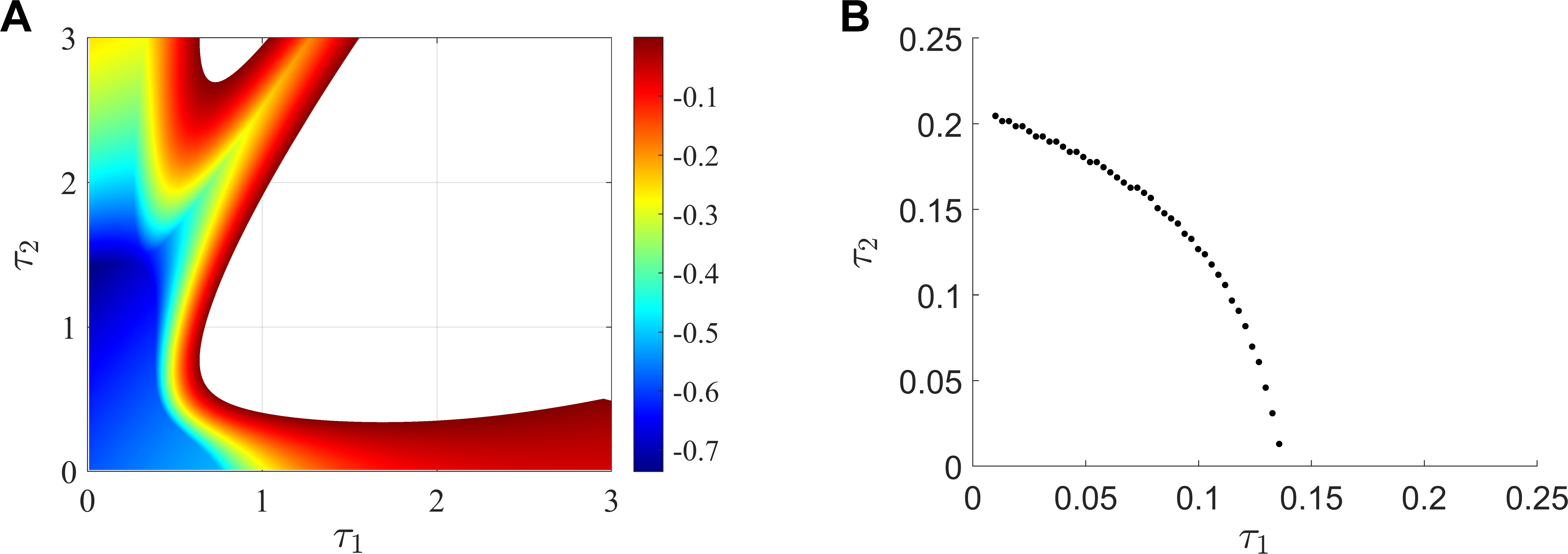}
    \end{center}
    \caption{Stability analysis for the second-order DDE with two delays (Eq.~\eqref{SO_2delaywd_eq1}) using parameter set 3: (A) stability chart obtained using the CCR method and (B) all points at which the Jacobian becomes rank-deficient.}
    \label{ex4fig1}
\end{figure}

The stability chart generated for parameter set 3 using the CCR method is shown in Fig.~\ref{ex4fig1}(A).
All points at which the Jacobian becomes rank-deficient (i.e., where $\left| \boldsymbol{J} \right| = 0$) while determining the stability regions are shown in Fig.~\ref{ex4fig1}(B).
At each of these 43 points, the ODEs become discontinuous and a new integration process is initiated at the following grid point with a fresh set of initial conditions, determined using the Galerkin approach.
This procedure adequately addresses the issue of encountering non-invertible Jacobians during continuation, and enables accurate determination of the stability regions despite these discontinuities.

\subsection{Computation time}\label{timecomp}
Finally, we report the computation time required to generate the stability charts shown in Figs.~\ref{ex1fig1}--\ref{ex3fig2} using the CCR method and the Galerkin approach.
All simulations were performed using \textsc{Matlab} R2018b on a 2.6-GHz Intel Xeon E5-2670 processor with 48 Gb of memory.
As shown in Table~\ref{Table1}, the stability charts presented here were generated between 3.9 and 10.3 times faster using the CCR method.
Furthermore, in the Galerkin approach, one must solve an eigenvalue problem of size $N_G \times N_G$ for a first-order system and $2N_G \times 2N_G$ for a second-order system at each grid point in the parametric space.
Note that the grid size has a significant effect on the simulation time in the Galerkin approach: simulation time increases linearly with the number of grid points.
Indeed, for any strategy in which an eigenvalue problem is solved at each grid point, the anticipated computation time is approximately $nmC$ for a grid of size $n \times m$, where $C$ is the computational cost of solving each eigenvalue problem.
In contrast, the CCR method is a continuation technique and the characteristic roots are obtained by solving a system of linearly independent ODEs, requiring substantially less computational effort.
As illustrated in Algorithm~1, the complexity of the CCR method is $N N_{p1} O(N_{p2})$, where $N$ is the number of roots being continued, $N_{p1}$ is the number of grid points over the domain of parameter 1, and $O(N_{p2})$ is the computational cost of solving a system of ODEs (Eq.~\eqref{mm_eq6}) over the domain of parameter 2.
The grid size does not dramatically affect the computation time in the CCR method and, as a result, stability regions can be readily determined with very high accuracy.
\begin{table}[htp]
    \caption{Computation time required to generate stability charts for DDEs.}
    \label{Table1}
    \centering
    \begin{tabular}{|l|c|c|}
    \hline
    \multirow{2}{*}{System} & \multicolumn{2}{c|}{Computation time (s)}\\
    \cline{2-3}
    & Galerkin approach & CCR method\\
    \hline
    Figure~\ref{ex1fig1} (Eq.~\eqref{fo_eq1}) & 688 & 95 (7.2$\times$ faster)\\
    \hline
    Figure~\ref{ex2fig1} (Eq.~\eqref{SO_eq1}) & 1083 & 105 (10.3$\times$ faster)\\
    \hline
    Figure~\ref{ex3fig1} (Eq.~\eqref{SO_2delaywd_eq1}, parameter set 1) & 653 & 108 (6.0$\times$ faster)\\
    \hline
    Figure~\ref{ex3fig2} (Eq.~\eqref{SO_2delaywd_eq1}, parameter set 2) & 662 & 168 (3.9$\times$ faster)\\
    \hline
    \end{tabular}
\end{table}

\section{Conclusions}\label{sec.conclusions}
We have developed a continuation of characteristic roots (CCR) method to determine the roots of the characteristic equation and obtain highly accurate stability charts for delay differential equations (DDEs) with multiple delays.
In this method, we write the characteristic roots as implicit functions of the parameter of interest and derive a continuation equation in the form of an ordinary differential equation (ODE).
The roots of the characteristic equation are determined by numerically integrating this derived system of linearly independent ODEs, using the solution obtained from the Galerkin approximation method as initial conditions.
The stability regions of the DDE are then determined by identifying the location of the rightmost characteristic root over the entire parametric space.
A key advantage of the proposed CCR method is that, rather than evaluating a large number of eigenvalue problems, highly accurate stability charts of the DDE are determined by solving a system of linearly independent ODEs.
Furthermore, the CCR method reduces the required computational time by a significant amount when compared to other available methods.
The efficacy of the proposed method has been demonstrated using first- and second-order DDEs with multiple delays.
The stability charts obtained in this work using the CCR method match those obtained using the Galerkin approach, and were generated between 3.9 and 10.3 times faster using the CCR method.
Finally, we have identified a limitation of the CCR method and recommended a technique to overcome rank-deficient Jacobians.
Although we limited our analysis to first- and second-order DDEs in this work, the CCR method can also be applied to generate stability charts for higher-order DDEs.


\section*{Acknowledgements}
Funding was provided to C.P.V.\ by the Department of Technology through the Inspire fellowship (grant number DST/INSPIRE/04/2014/000972).
The funders had no role in study design, data collection and analysis, decision to publish, or preparation of the manuscript.

\section*{Competing Interests}
The authors declare that no competing interests exist.

\end{document}